\documentclass{article}
\usepackage{amsthm,amsfonts, amsbsy, amssymb,amsmath,graphicx}
\usepackage{graphics}

\newtheorem{thm}{Theorem}

\newtheorem{crl}{Corollary}

\title{Free Knots and Groups}

\author{V.O.Manturov \footnote{partially supported by RFBR No. 07-01-00648}, O.V.Manturov}
\date{}

\begin{document}

\maketitle

Virtual knot theory invented by Kauffman \cite{KaV} is an important
generalization of knot theory; some methods of classical knot theory
can be generalized to virtual knot theory straightforwardly, and
some other can not, \cite{Mybook}. On the other hand, there are lots
of new constructions coming from virtual knot theory and similar
theories, e.g., graph-links by Ilyutko and V.O.Manturov \cite{IM}.
In the present paper we consider {\em free knots}, a thorough
simplification of virtual knots, first introduced by Turaev
\cite{Tu} under the name of {\em homotopy classes of Gauss words}.
Turaev conjectured the non-triviality of free knots, and the first
examples of non-trivial free knots were constructed by the first
named author, \cite{Ma} and A.Gibson \cite{Gib}. In \cite{Ma},
several important theorems about free knots were proved by using the
notion of {\em parity}: a chord in the chord diagram is {\em even}
if the number of chords it is linked with, is even; otherwise it is
{\em odd}. Free knots are equivalence classes of Gauss diagrams
(chord diagrams) by the relations corresponding to the three
Reidemeister moves.

 Non-triviality of free knots yields non-triviality of the
underlying virtual knots and usually allows to improve many of
virtual knot invariants (since the discovery of parity, several new
applications appeared in \cite{Af,CM,MTrst}.

Recently, a partial case of the invariant constructed in the present
paper was proved to be a {\em sliceness obstruction} for free knots,
see \cite{Sl}.

 In the
present paper we construct a simple and rather  strong invariant of
free knots, valued in a certain group (more precisely, there will be
a group for every natural $m$, and the invariants for groups with
greater $m$ naturally generalize invariants for smaller $m$). This
invariant is constructed only out of the notion of parity, and in
fact this invariant depends merely on the disposition of the chord
ends rather than on the chord diagram itself.

 Let $C$  be the segment  of positive integers
  from $1$ to a given integer $2n$. An unordered partition $a=
\{(p_1,q_1),(p_2,q_2)...(p_n,q_n)\}$ of $C$   (all numbers
 $p_i,q_j,i,j=1,2,...,n$ are distinct) is called  a
 \it{chord diagram (with a base point).}
 \rm Each pair in $a$ is  a \it{chord}, \rm each element of pair is
 called an end of chord. We say that two chords $(p_i,q_i),(p_j,q_j)$ are
  {\em linked} (resp., {\em unlinked})  depending on whether the  following
 statement is true or false:
$$(p_i-p_j)(p_i-q_j)(q_i-p_j)(q_i-q_j)<0.$$

Here two chords $(p_{1},p_{2})$ and $(q_{1},q_{2})$ are {\em linked}
whenever two half-circles connecting $(p_{1},0)$ to $(p_{2},0)$ and
$(q_{1},0)$ to $(q_{2},0)$ in the upper half-plane have an
intersection point. Note that the property of being linked does not
change under the cyclic permutation of partitioned chord ends $1\to
2\to  \dots \to 2n\to 1$.

Fix a positive integer $m$.

Let $g(p,b)$ be the number of all chords belonging to given set of
chords $b$ and linked with the chord $p$. We denote by $a_0$ the
following subset of the set $a$:
            $$a_0=\{p\in a|g(p,a)\ is \ odd\}$$
We proceed by induction (for $k\le m-1$) to define  $a_k$ as
$$a_k=\{p \in a\backslash \cup_{i<k}a_i|g(p,a\backslash
\cup_{i<k}a_i) \ is\ odd\}$$
 Here $a\backslash \cup_{i<k}a_i$ denotes the complement to
 $\cup_{i<k}a_i$ in $a$.

Finally, $a_{m}$ denotes the complement to $\cup_{j=1}^{m-1} a_{j}$
in $a$: these are all chords which survive after $m$ consecutive
operations of deleting odd chords.

 We split each set $a_k, k=0,1,\dots ,m-1$  into two
 disjoint subsets $a_k=a'_k \sqcup a''_k $,  by putting
 $$a'_k= \{p \in a_k|g(p,a_k) \ is \ odd \}$$. Having a chord diagram $c$, we construct
a word $w(c)$ in the alphabet
$$D=\{a'_0,a''_0,a'_1,a''_1,...,a'_{m-1},a''_{m-1},a_m\}$$
as follows. The letter number $k$ in the word $w(c)$ to be
constructed will be denoted by the same letter ($a'_{j}, a''_{j}$ or
$a_j$ if $j=n$) as the subset of chords, the corresponding chord
belongs to. We say that the end of a chord is marked with letter $j$
from the alphabet $D$.

Each word in  $D$ can be considered as an element of some group $G$
generated by $D$. Define the group $G$ by  generators from $D$ and
the following relations
 $$a'_0 a'_0=e,a''_0 a''_0=e,a'_1 a'_1=e, a''_1 a''_1=e,\dots,
a'_{m-1} a'_{m-1}=e, a''_{m-1} a''_{m-1}=e,a_{m}a_{m}=e,$$
 $$ a''_i a'_j=a'_j a'_i,i<j,
a''_i a''_j=a''_j a'_i, i<j, a'_{i}a_{m}=a_{m}a''_{i}, i<m.
\eqno(1)$$ Here $e$ denotes the unity element in $G$.

Here we abuse the notation by omitting the dependence of $G$ on $m$;
here $m$ is fixed once forever.

 Now we are going to define free knots as equivalence classes of chord diagrams
 modulo ``Reidemeister moves''
on the set of all chord diagrams. The main statement of our work is

\begin{thm}
If two chord diagrams $c_1,c_2$ are equivalent (represent the same
long free knot) then we have
$$w(c_1)=w(c_2)$$ in $G$.
\end{thm}

A {\em long free knot} is an equivalence class of chord diagrams
with a base point by {\em Reidemeister moves}. These Reidemeister
moves correspond to usual Reidemeister moves applied to Gauss
diagrams if we forget the information about arrows and signs
corresponding to chords.

The first increasing (resp., decreasing) Reidemeister move is an
addition (resp.,  removal) to given chord diagram $c$ of a chord
$(p,q)$ such that $|p-q|=1$
 The second increasing (resp., decreasing) Reidemeister move is addition/ removal a pair of
``adjacent''\; chords. Two chords $(p_i,q_i),(p_j,q_j)$  of diagram
$c$ are  {\em adjacent} if
$|min(p_{i},q_{i})-min(p_{j},q_{j})|=|max(p_{i},q_{i})-max(p_{j},q_{j})|=1$
 (in both cases, for the first and the second Reidemeister moves, after such an addition/removal, the
remaining chords are renumbered accordingly).

We say that a triple of chords
$(p_{i},q_{i}),(p_{j},q_{j}),(p_{k},q_{k})$ is {\em completely
adjoint} if the six ends of these chords can be partitioned such a
way that each pair contains two ends of different lower indices
$i,j,k$ and the two numbers in the each pair differ by one, e.g.,
$|p_{i}-p_{j}|=|p_{k}-q_{i}|=|q_{j}-q_{k}|=1$.

 The third Reidemeister move is defined only for those diagrams, which
contains a completely adjoint   triple. It is easy to see that the
six ends of a completely adjoint triple represent an set of integers
looking like
$$S=\{r,r+1,s,s+1,t,t+1\}, \eqno(2)$$ and elements $r,r+1$ (and $s,s+1$  and
$t,t+1$) belong to different chords. We define an involution
$f:S\rightarrow S$ by setting
$$r\leftrightarrow {r+1},s\leftrightarrow {s+1}, t\leftrightarrow {t+1}\eqno(3)$$

 Now we define the transformation of a triple $T$ into  triple $T'$
according the following rule: the set of chord ends of the triple
$T$ coincides with that of triple $T'$, and a pair of integers $u,v
\in S$ forms  a chord in $T'$ whenever the pair $f(u),f(v)$ forms a
chord in $T$.  The third
 Reidemeister move transforms a given  chord diagram including a completely adjoint
 triple $T$  into the diagram obtained by replacing the triple $T$ by $T'$ and leaving
 the remaining chords fixed.

So, we have completed the definition of a long free knot.

To prove Theorem 1, we check the invariance of $w(c)$ under
Reidemeister's moves.

For the first Reidemeister move $c_{1}\to c_{2}$, $w(c_1)$  and
$w(c_2)$ are obtained from each other by an addition/removal of a
couple of consequitve identical letters (generators of the group),
which yields the identity in the group.

In the case of the second  Reidemeister move $c_{1}\to c_{2}$, we
add two pairs of identical letters, namely, follows that the words
corresponding to diagrams $c_1,c_2$ look like
 $$w(c_1)=UVW;  w(c_2)=Uz_1z_2Vz_3z_4W,\eqno(4)$$
 where $U,V,W$ are some subwords form the word
 $w(c_1)$, and $z_1,z_2,z_3,z_4$ are
 letters from  $D$ corresponding to the ends of  $u,v$.
We have $z_1=z_2=z_3=z_4$, so $w(c_1)=w(c_2)$ in $G$.

For the third Reidemeister move, we get a word of the following type
$$w(c_1)=U\alpha V\beta W\gamma X, \eqno(5)$$
constructed according to the above  rules for the chord diagram
$c_1$ and the word
$$w(c_2)=U\delta V\epsilon W \zeta X, \eqno(6)$$
constructed from the  diagram $c_2$; the latter is obtained from
$c_1$ by means of the third Reidemeister move.

Here each of $\alpha,\beta,\gamma,\delta,\epsilon,\zeta$ is a pair
of generators of the group $G$ corresponding to the two adjacent
chord ends in an adjoint triple.

Our goal is to show that in $G$ the following equalities hold:
$\alpha=\delta$, $\beta=\epsilon$, $\gamma=\zeta$.

Every chord diagram $D$ containing a triple of completely adjoint
chords has the following property: the number of odd chords in the
triple is {\em even}, i.e., is equal to zero or two. This follows
from the Pasch axiom of the Hilbert axiom system. (A chord $b\in c$
is {\em even (odd)} iff $g(b,c)\equiv 0\; mod\; 2$ (resp., $\equiv\
1 \; mod \; 2$)). An analogous property takes place for every of set
from list $a_k,k=0,1,\cdots, m-1$. This means: if any of
$a_k,k=0,1,\cdots, m-1$ contains a completely adjoint triple, then
the number of chords $b$ in the completely adjoint triple,
satisfying the condition $g(b,a_k)\equiv 1 \ mod \; 2,$ is even.

Thus, we see that the triple of completely adjoint chords contains
either two odd chords or zero odd chords.

Denote by $h(c)$ the chord diagram obtained from the diagram $c$ by
deleting all odd chords in it. If $c_1$ has a completely adjoint
triple without odd chords in it, then we pass to $h(c_1)$. The
deletion of all odd chords leaves the property of a triple (of
persistent chords) to be completely adjoint true. If  $h(c_1)$ has a
completely adjoint triple without odd chords in it then we pass to
$h(h(c_1))$, and so on. As a result we get the following two
options:

1. We obtain a diagram $c^{*}_1$ with exactly two odd chords in the
triple.

2. We obtain a diagram $c^{**}_2$ without odd chords at all.

In the second case, all chords of the initial completely adjoint
triple in $c_{1}$ have index $m$. So are the corresponding chords
from $c_{2}$. So, each of the words
$\alpha,\beta,\gamma,\delta,\epsilon,\zeta$ looks like $a_{m}\cdot
a_{m}$, and the claim follows. So, the words
 $w(c_1)$ and $w(c_2)$ identically coincide.

 It remains to consider the first case.
In this case the sequence $c_{1},h(c_{1}),h(h(c_{1})),\dots$
contains a chord diagram $c^{*}_1$ including a completely adjoint
triple $T$ with exactly two odd chords. The six letters which mark
the ends are coupled into elements $\alpha,\beta,\gamma\in G$.
Without loss of generality, assume $\alpha$ is the the product of a
pair of generators corresponding to two adjacent ends of odd chords.
So, $\alpha$ is equal to one of the following: $a'_k\cdot a'_k$ or
$a'_k \cdot a^{''}_k$ or $a^{''}_k \cdot a'_k$ or $a^{''}_k \cdot
a^{''}_k$.

Notice that the indices of chords in $\delta$ will be the same as in
$\alpha$. Now, it is easy to see that either the corresponding
subword $\delta$ in $w(c_2)$ will be the same if $\alpha$ is equal
to one of $a'_{k}\cdot a''_{k}$ or $a''_{k}\cdot a'_{k}$. In case
when $\alpha$ is a square of a generator $a'_{k}$ (resp.,
$a''_{k}$), the word $\delta$ is a square of the other generator
$a''_{k}$ (resp., $a'_{k}$). So, $\alpha=\delta$ in $G$.

Now consider another segment of $w(c_1)$ (say, $\beta$) and the
corresponding segment of $w(c_2)$ (in this case $\epsilon$).

In $\beta$ one letter is  $a'_{k}$ or $a''_{k}$ and the other chord
end with a higher index, then when passing from $\beta$ to
$\epsilon$, the generators change their places and the generator
with smaller index transforms $a'_{k}\longleftrightarrow a''_{k}$,
whence the generator with higher index remains the same.

So, the equality $\beta=\epsilon$ is one of the relations of $G$.
Analogously, $\gamma=\zeta$.

The theorem is proved.

 Chord diagrams considered above deal with
so-called  long free knots, i.e. free knots with a chosen initial
point.

{\em Free knots} are equivalence classes of long free knots by the
move which {\em changes the initial point}; this move acts by a
cyclic permutation $1\to 2\to \dots \to 2n\to 1\to 2n$ on the set of
partitioned points.

It is obvious that if $c_1$ is obtained from $c_2$ by such an
operation then $w(c_1)$ and $w(c_2)$ are conjugate in $G$. This
yields the following

\begin{crl}
The conjugacy class $[w(\cdot)]$ in $G$ is the invariant of free
knots.
\end{crl}

For a given $m$, the group $G$ has a very simple Cayley graph.
Namely, elements of $G$ are in one-to-one correspondence with the
set of points in Euclidean space ${\bf R}^{m+1}$ such that their
coordinates are integers and the last coordinate is $0$ or $1$.

Here, the origin of coordinates corresponds to the unit of the group
$G$. The right muliplication by element with low index $k$ (i.e.,
$a'_{k}$ or $a''_{k}$ or $a_{k}$ if $k=m$) corresponds to one step
shift along the coordinate number $k$ defined as follows. Let
$x_1,x_2,...,x_{m+1}$ be the coordinates of given point.

 The multiplication by $a'_0$ on the right ($a''_0$) increases (decreases) the
  first coordinate $x_1$ if $\sum_{s=1}^m
x_s$ is even (odd), the multiplication by  $a'_1$ ($a''_1$)
increases (decreases) the second coordinate $x_2$ if $\sum_{s=2}^m
x_s$ is even (odd); the multiplication by $a'_k$ ($a''_k$) for $k\le
m$ increases (decreases) the coordinate number  $x_{k+1}$ if
$\sum_{s=k}^m x_s$ is even (odd) and. Finally,  $a_m$ changes the
coordinate $x_{m+1}$ from zero to one and from one to zero.

\begin{figure}
\centering\includegraphics[width=200pt]{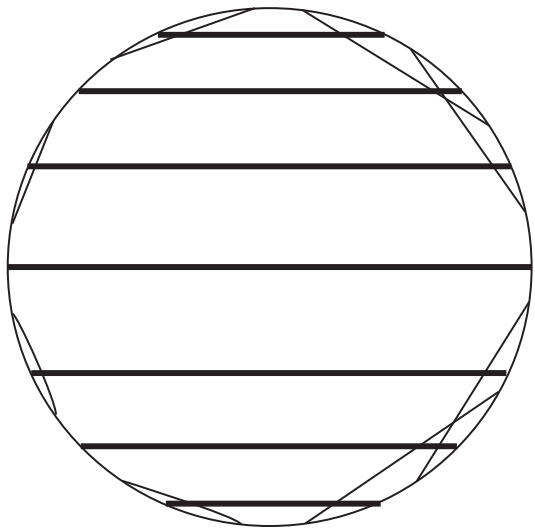} \caption{A
non-trivial free knot} \label{fig1}
\end{figure}

In Fig. \ref{fig1}, we present a non-trivial free knot recognizable
by the group with $m=1$.

  \end{document}